\documentclass[11pt]{article}

\usepackage[T2A]{fontenc}
\usepackage[cp1251]{inputenc}
\usepackage{amsthm}
\usepackage{amsfonts}
\usepackage{eufrak}
\usepackage{amssymb}
\usepackage{amsmath}
\usepackage{cite}

\makeatletter
\renewcommand{\@biblabel}[1]{#1.\hfill}
\makeatother

\begin{document}
\newtheoremstyle{mytheorem}
  {\topsep}   
  {\topsep}   
  {\itshape}  
  {}       
  {\bfseries} 
  { }         
  {5pt plus 1pt minus 1pt} 
  { }          
\newtheoremstyle{myremark}
  {\topsep}   
  {\topsep}   
  {\upshape}  
  {}       
  {\bfseries} 
  {  }         
  {5pt plus 0pt minus 1pt} 
  { }          
\theoremstyle{mytheorem}
\newtheorem{theorem}{Theorem}[section]
 \newtheorem{theorema}{Theorem}
  \newtheorem*{heyde*}{The Heyde theorem}
 \newtheorem*{a*}{The Skitovich--Darmois theorem}
 \newtheorem*{b*}{The Skitovich--Darmois theorem for the group $\mathbb{R}\times \mathbb{Z}(2)$}
\newtheorem{proposition}{Proposition}[section]
\newtheorem{lemma}{Lemma}[section]
\newtheorem{corollary}{Corollary}[section]
\newtheorem{definition}{Definition}[section]
\theoremstyle{myremark}
\newtheorem{remark}{Remark}[section]
\noindent This article is accepted for publishing in

\noindent the journal Potential Analysis 

\vskip 1 cm

\noindent{\bf Heyde theorem on locally compact    Abelian groups }

\noindent{\bf  with  the connected component of zero of dimension 1}

\bigskip

\noindent{\textbf{Gennadiy Feldman}}

\bigskip

\noindent{\textbf{Abstract}}

\bigskip

\noindent  Let $X$ be a locally compact    Abelian group with the connected component of zero of dimension 1. Let $\xi_1$ and $\xi_2$ be independent random variables with   values in   $X$ with nonvanishing characteristic functions. We prove that if  a topological automorphism $\alpha$ of the group $X$ satisfies the condition
${{\rm Ker}(I+\alpha)=\{0\}}$
and the conditional  distribution of the linear form ${L_2 = \xi_1 + \alpha\xi_2}$ given ${L_1 = \xi_1 + \xi_2}$ is symmetric, then the distributions of  $\xi_j$  are   convolutions of Gaussian distributions on  $X$ and distributions supported in the subgroup $\{x\in X:2x=0\}$. This result can be viewed as a generalization of the well-known Heyde theorem on the characterization of the Gaussian distribution on the real line.

\bigskip


\bigskip


\section{ Introduction}

\noindent According to the well-known Heyde  theorem the Gaussian distributions on the real line are characterized by the symmetry of the conditional distribution of one linear form of independent random variables  given another (see \cite[\S\,13.4.1]{KaLiRa}).  A number of papers is devoted to group analogues  of this theorem in the case when
  independent random variables
   take values in various locally compact
   Abelian groups,
and   coefficients of linear forms are
topological automorphisms of a group  (see
  \cite{Fe2, Fe4,  Fe6, Fe20bb, M2013, M2020, POTA, F_solenoid,  JFAA2021,  FeTVP1, FeTVP, My2,   MiFe1}, and also   \cite[Chapter VI]{Fe5}). It   turns out that the class of distributions, which is characterized by the symmetry of the conditional distribution of one linear form given another, depends on the presence of elements of order 2 in the group and conditions on  coefficients of the linear forms and on the distributions.  In particular, it is   important whether the characteristic functions of the considered distributions vanish or not.
In this article, we continue these investigations.

Let $X$ be a locally compact    Abelian group with the connected component of zero of dimension 1, and let $\xi_1$ and $\xi_2$  be independent random variables  with   values in the group   $X$  with nonvanishing characteristic functions. Denote by $I$ the identity automorphism of a group.
The main result of the article is the following.  If  $\alpha$ is a topological automorphism of the group  $X$ such that
\begin{equation}\label{1}
{\rm Ker}(I+\alpha)=\{0\},
\end{equation}
and the conditional  distribution of the linear form $L_2 = \xi_1 + \alpha\xi_2$ given $L_1 = \xi_1 + \xi_2$ is symmetric, then the distributions of  $\xi_j$ are convolutions of Gaussian distributions on  $X$ and distributions supported in the subgroup generated by all elements of order 2 of the group $X$.
It is interesting to note that the proof of this result uses some properties of entire functions, whereas the proof of the Heyde theorem on the real line for two independent random variables does not use   complex variables. Note also that the proof of the main theorem is reduced to a description of the solutions of some functional equation on the character group of the group $X$.

We  need some notions and results from the duality theory for locally compact Abelian groups   (see \cite{Hewitt-Ross}). Let $X$ be a second countable locally  compact Abelian group.  We will consider only such groups. Denote by $Y$ the character
group of the group $X$, and by $(x,y)$ the value of a character $y \in Y$ at an element  $x\in X$. If $H$ is a  subgroup of the group  $Y$, denote by   $A(X, H) = \{x \in X: (x, y) = 1$ for all $y \in H \}$
its annihilator. Let $\alpha$ be a topological automorphism of the group $X$.  If $G$ is a closed subgroup of the group $X$  and
 $\alpha(G)=G$, then the restriction of $\alpha$  to $G$ is a topological automorphism of the group $G$. Denote by $\alpha_G$ this automorphism.
Let  $n$ be an integer, $n\ne 0$. Denote by $f_n:X \rightarrow X$ an  endomorphism
of the group $X$
 defined by the formula  $f_nx=nx$, $x\in X$. Put $X^{(n)}=f_n(X)$. 
 Let $X_1$ and $X_2$ be locally compact Abelian groups with the character groups $Y_1$ and $Y_2$ respectively.
If $\alpha:X_1\rightarrow X_2$ is a continuous homomorphism,  then the adjoint homomorphism  $\tilde\alpha:Y_2\rightarrow Y_1$  is defined by the formula
 $(x_1, \tilde\alpha y_2) = (\alpha x_1, y_2)$ for all $x_1 \in X_1$, $y_2 \in Y_2$. Denote the connected component of zero of the group $X$ by $c_X$.
Denote by $\mathbb{R}$ the additive group of real numbers, by  $\mathbb{Z}(2)=\{0, 1\}$ the group of residue classes modulo 2, by  $\mathbb{T}$ the circle group, by $\mathbb{Z}$ the set of integers and by $\mathbb{C}$ the complex plane.

Let $f(y)$ be an arbitrary function on the group $Y$, and let $h \in
Y$.  Denote by $\Delta_h$ the finite difference operator
$$
\Delta_h f(y) = f(y + h) - f(y), \ \ y \in Y.
$$
A function $f(y)$ on the group $Y$ is called a polynomial if for a natural $n$ the function $f(y)$ satisfies the equation
$$
\Delta_h^{n} f(y) = 0, \ \ y, h \in Y.
$$

Denote by ${\rm M}^1(X)$  the convolution semigroup of probability distributions on the group   $X$. Let  $\mu\in{\rm M}^1(X)$. Denote by $$\hat\mu(y) =
\int_{X}(x, y)d \mu(x), \ \ y\in Y,$$
the characteristic function (Fourier transform) of
the distribution  $\mu$, and by $\sigma(\mu)$ the support of $\mu$. Define the distribution $\bar \mu \in {\rm M}^1(X)$ by the formula
 $\bar \mu(B) = \mu(-B)$ for any Borel  subset $B$ of $X$.
Then $\hat{\bar{\mu}}(y)=\overline{\hat\mu(y)}$. For $x\in X$ denote by $E_x$  the degenerate distribution
 concentrated at the element $x$.

 A distribution  $\gamma$ on the   group  $X$ is said to be Gaussian
(\!\!\cite[Chapter IV]{Pa})  if its characteristic function is represented in the form
$$
\hat\gamma(y)= (x,y)\exp\{-\varphi(y)\}, \ \  y\in Y,
$$
where $x \in X$, and $\varphi(y)$ is a continuous nonnegative function on the group   $Y$ satisfying the equation
$$
\varphi(u + v) + \varphi(u
- v) = 2[\varphi(u) + \varphi(v)], \ \ u,  v \in
Y.
$$
Observe that, according to this definition, the degenerate distributions are Gaussian. Denote by  $\Gamma(X)$ the set of Gaussian distributions on the group  $X$.

\section{Heyde  theorem for the group $\mathbb{R}\times D$,  where $D$ is a locally compact  totally disconnected Abelian group}

 Let $X$ be a locally compact Abelian group, and let $\alpha_j, \beta_j$ be topological automorphisms of   $X$. Let
$\xi_1$ and $\xi_2$ be independent random variables with values  in the group  $X$ and distributions $\mu_1$ and $\mu_2$. Assume that the conditional distribution of the linear form $L_2 = \beta_1\xi_1 + \beta_2\xi_2$ given $L_1 = \alpha_1\xi_1 + \alpha_2\xi_2$ is symmetric. It is easy to see that the description of possible distributions  $\mu_j$ is reduced to the case when     $L_1 = \xi_1 + \xi_2$ and $L_2 = \xi_1 + \alpha\xi_2$, where $\alpha$ is a topological automorphism  of the group $X$.

At the beginning we will prove an analogue of the Heyde theorem for the group of the form $\mathbb{R}\times D$,  where $D$ is a locally compact  totally disconnected Abelian group. The proof of the main theorem is based on it.

\begin{theorem}\label{th1} Let $X$ be a group of the form   $X=\mathbb{R}\times D$, where  $D$ is a   locally compact  totally disconnected Abelian group. Set $G=\{x\in D:2x=0\}$. Let  $\alpha$ be a topological automorphism of the group    $X$ satisfying condition $(\ref{1})$. Let $\xi_1$ and $\xi_2$ be independent random variables with   values in   $X$ and distributions
$\mu_1$ and $\mu_2$ with nonvanishing characteristic functions.
If the conditional distribution of the linear form    $L_2 = \xi_1 + \alpha\xi_2$ given $L_1 = \xi_1 + \xi_2$ is symmetric, then
$\mu_j=\gamma_j*\omega_j*E_{x_j}$, where $\gamma_j \in \Gamma({\mathbb R})$,
 $\omega_j\in {\rm M}^1(G)$,
$x_j\in X$, $j=1, 2$.
\end{theorem}

To prove Theorem \ref{th1}, we need the following lemmas.

\begin{lemma}\label{lem1}{\rm(\!\!\cite[Lemma 16.1]{Fe5})}  Let $X$ be a locally compact Abelian group, and let  $\alpha$ be a topological automorphism of  the group $X$. Let $\xi_1$ and $\xi_2$ be independent random variables with   values in   $X$ and distributions
$\mu_1$ and $\mu_2$. The conditional distribution of the linear form    $L_2 = \xi_1 + \alpha\xi_2$ given $L_1 = \xi_1 + \xi_2$ is symmetric
  if and only if the characteristic functions  $\hat\mu_j(y)$  satisfy the equation
\begin{equation}\label{2}
\hat\mu_1(u+v )\hat\mu_2(u+\tilde\alpha v )=
\hat\mu_1(u-v )\hat\mu_2(u-\tilde\alpha v), \ \ u, v \in Y.
\end{equation}
\end{lemma}
Lemma \ref{lem1} means that the description of possible distributions $\mu_j$ under the conditions of the Heyde theorem on the group $X$ is reduced to the description of solutions of functional equation (\ref{2}) in the class of characteristic functions on the group $Y$, that is, to some problem of abstract harmonic analysis.
\begin{lemma}\label{lem3} {\rm(\!\!\cite[Lemma 3.8]{M2013})}    Let $X$ be a locally compact Abelian group, and let $\alpha$ be a topological automorphism of  the group $X$. Let
$\mu_1$ and $\mu_2$ be distributions on  $X$ such that the characteristic functions $\hat\mu_j(y)$ satisfy equation  $(\ref{2})$. Then $\hat\mu_j(y)$ satisfy the equation
\begin{equation}\label{5}
\hat\mu_1((I+\tilde\alpha) u+2 v)\hat\mu_2(2\tilde\alpha  u+(I+\tilde\alpha) v)$$$$=\hat\mu_1((I+\tilde\alpha)
u)\hat\mu_2(2\tilde\alpha
u)\hat\mu_1(2 v)\hat\mu_2((I+\tilde\alpha) v), \ \ u, v \in
Y.
\end{equation}
\end{lemma}
\begin{lemma}\label{lem4new}  {\rm(\!\!\cite[(24.41)]{Hewitt-Ross})}
Let $X_1$ and $X_2$ be locally compact Abelian groups with the character groups $Y_1$ and $Y_2$ respectively.
Let $\alpha:X_1\rightarrow X_2$ be a continuous homomorphism, and let $\tilde\alpha:Y_2\rightarrow Y_1$ be the adjoin homomorphism. The subgroup $\tilde\alpha(Y_2)$ is dense in $Y_1$ if and only if
$\alpha$ is a monomorphism.
\end{lemma}
The following statement is well known. For the proof see e.g.     \cite[Proposition 5.7]{Fe5}.
\begin{lemma}\label{lem4}
Let $Y$ be a locally compact Abelian group  such that all its elements  are compact. Let $P(y)$ be a continuous polynomial on   $Y$. Then $P(y)=const$.
\end{lemma}
\begin{lemma}\label{lem6}
Let   $X$ be a locally compact Abelian group,   let $H$ be a closed subgroup of the group $Y$, and let $\mu$ be a distribution on the group   $X$ such that $\hat\mu(y)=1$ for all $y\in H$.  Then $\sigma(\mu)\subset A(X, H)$.
\end{lemma}
\begin{lemma}\label{lem5} {\rm(\!\!\cite[Lemma 5]{Fe4})} Let $X$ be a group of the form
$X = {\mathbb R}  \times G$, where  $G$ is a locally compact Abelian group such that all its nonzero elements  have order
 $2$. Let
$\mu$ be a distribution on $X$ such that  $\mu =\gamma*\omega$, where $\gamma
\in \Gamma({\mathbb R})$, $\omega \in {\rm M}^1(G)$, and the characteristic function
of the distribution $\omega$ does not vanish. Let $\mu=\mu_1*\mu_2$, where $\mu_j$ are distributions on  $X$. Then   $\mu_j = \gamma_j*\omega_j$, where $\gamma_j \in \Gamma({\mathbb
R})$, $\omega_j\in {\rm M}^1(G)$, $j=1, 2$.
\end{lemma}
Note that if    $G$  is a locally compact
Abelian group such that every element of $G$ different from zero   has
order  $2$,   then  $G$ is topologically isomorphic the a group of the form
$$ {\mathbb{Z}(2)}^{\mathfrak{n}}\times {\mathbb{Z}(2)}^{\mathfrak{m}*},$$
where $\mathfrak{n}$ and $\mathfrak{m}$ are cardinal numbers, the group   ${\mathbb{Z}(2)}^{\mathfrak{n}}$ is considered in the product topology, and  the group ${\mathbb{Z}(2)}^{\mathfrak{m}*}$ is considered in the discrete topology  (\!\!\cite[$($25.29)]{Hewitt-Ross}).

We formulate as a lemma  the Heyde theorem on the real line for two independent random variables.
\begin{lemma}\label{lem7}   Let $\xi_1$ and $\xi_2$ be independent random variables with     distributions
$\mu_1$ and $\mu_2$.   Let
$a\ne 0$, $a\ne -1$. If the conditional distribution of the linear form    $L_2 = \xi_1 + a\xi_2$ given $L_1 = \xi_1 + \xi_2$ is symmetric, then
   $\mu_j\in\Gamma(\mathbb{R})$, $j=1, 2$.
\end{lemma}

Let $X$ be a group of the form
$X = {\mathbb R}  \times G$, where  $G$ is a locally compact Abelian group. Then the group $Y$ is of the form   $Y=\mathbb{R}\times H$, where the group     $H$   is topologically isomorphic to the character group of the group    $G$. Denote by $(s, h)$, where $s\in \mathbb{R}$, $h\in H$, elements of the group  $Y$. We   need the following easily verifiable assertion  (see e.g.   \cite[Lemma 6.9]{Fe9}).
\begin{lemma}\label{lem2}
Let $X$ be a group of the form
$X = {\mathbb R}  \times G$, where  $G$ is a locally compact Abelian group. Let
$\mu$ be a distribution on $X$. Assume that the function $\hat\mu(s, 0)$, $s\in \mathbb{R}$, can be extended  to the complex plane $\mathbb{C}$  as an entire function in  $s$. Then for each fixed      $h\in H$ the function $\hat\mu(s, h)$, $s\in \mathbb{R}$, can be also extended  to the complex plane $\mathbb{C}$ as an entire function in $s$, the function $\hat\mu(s, h)$ is continuous on $\mathbb{C}\times H$,   and the inequality
$$
\max_{s\in \mathbb{C}, \ |s|\le r}|\hat\mu(s, h)|\le \max_{s\in \mathbb{C}, \ |s|\le r}|\hat\mu(s, 0)|, \ \ h\in H,
$$
holds.
\end{lemma}
\noindent{\it Proof of Theorem \ref{th1}}  We divide the  proof of the theorem into several steps.

1. Put $\nu_j = \mu_j* \bar \mu_j$, and first prove that
\begin{equation}\label{18.04.4}
\sigma(\nu_j)\subset \mathbb{R}\times G, \ \ j=1, 2.
\end{equation}
It is obvious that the group $Y$ is of the form   $Y=\mathbb{R}\times K$, where the group     $K$   is topologically isomorphic to the character group of the group    $D$.
 By Lemma \ref{lem1}, the symmetry of the conditional distribution of the linear form
 $L_2$ given
  $L_1$ implies that the characteristic functions
$\hat\mu_j(y)$ satisfy equation   $(\ref{2})$. Obviously,   $\hat \nu_j(y) = |\hat \mu_j(y)|^2>0$ for all  $y \in Y$, $j=1, 2$, and
   the characteristic functions
 $\hat \nu_j(y)$
also satisfy equation   $(\ref{2})$.
By Lemma  \ref{lem3}, it follows from this that  the characteristic functions  $\hat\nu_j(y)$ satisfy equation $(\ref{5})$.
Put $\psi_j(y)=-\log\hat\nu_j(y)$, $j=1, 2$, and note that
$(\ref{5})$ implies that the functions  $\psi_j(y)$  satisfy the equation
\begin{equation}\label{6}
\psi_1((I+\tilde\alpha) u+2 v)+\psi_2(2\tilde\alpha  u+(I+\tilde\alpha) v)$$$$=\psi_1((I+\tilde\alpha)
u)+\psi_2(2\tilde\alpha
u)+\psi_1(2 v)+\psi_2((I+\tilde\alpha) v), \ \ u, v \in
Y.
\end{equation}
Put
\begin{equation}\label{16.09.15.1}
P(y)=\psi_1((I+\tilde\alpha)
y)+\psi_2(2\tilde\alpha
y), \ \ Q(y)=\psi_1(2 y)+\psi_2((I+\tilde\alpha) y), \ \ y\in Y,
\end{equation}
and rewrite equation (\ref{6}) in the form
\begin{equation}\label{09.02.2}
\psi_1((I+\tilde\alpha) u+2 v)+\psi_2(2\tilde\alpha  u+(I+\tilde\alpha) v)=P(u)+Q(v), \ \ u, v \in
Y.
\end{equation}
To solve  equation  (\ref{09.02.2}) we use the finite-difference method (see e.g. \cite[\S 10]{Fe5}). This is the standard reasoning, and we will present it for completeness.

Take an arbitrary element $h_1$ of the group $Y$.
Substitute   $u+(I+\tilde\alpha) h_1$ instead of $u$ and $v-2\tilde\alpha  h_1$ instead of $v$ in
 (\ref{09.02.2}). Subtracting
 equation
   (\ref{09.02.2}) from the resulting equation, we get
  \begin{equation}\label{7}
    \Delta_{(I-\tilde\alpha)^2 h_1}{\psi_1((I+\tilde\alpha) u+2 v)}
    =\Delta_{(I+\tilde\alpha) h_1} P(u)+\Delta_{-2\tilde\alpha  h_1} Q(v),
\ \ u, v\in Y.
\end{equation}
  Take an arbitrary element $h_2$ of the group $Y$.
Substitute   $u+2h_{2}$ instead of $u$ and $v-(I+\tilde\alpha)h_{2}$ instead of $v$ in
(\ref{7}). Subtracting
    equation
  (\ref{7}) from the obtaining equation, we have
 \begin{equation}\label{8}
     \Delta_{2 h_2}\Delta_{(I+\tilde\alpha) h_1} P(u)+\Delta_{-(I+\tilde\alpha) h_2}\Delta_{-2\tilde\alpha  h_1} Q(v)=0,
\ \ u, v\in Y.
\end{equation}
Take an arbitrary element $h_3$ of the group $Y$.
Substitute   $u+h_3$ instead of $u$ in
(\ref{8}). Subtracting
    equation
  (\ref{8}) from the resulting equation, we get
 \begin{equation}\label{9}
   \Delta_{h_3}\Delta_{2 h_2}\Delta_{(I+\tilde\alpha) h_1} P(u)=0,
\ \ u\in Y.
\end{equation}
By Lemma \ref{lem4new}, it follows from  (\ref{1}) that the subgroup   $(I+\tilde\alpha)(Y)$ is dense in $Y$. Taking into account that   $h_j$ and $u$ in  (\ref{9})  are arbitrary elements of $Y$, from  (\ref{9})   we obtain that the function $P(y)$ satisfies the equation
\begin{equation}\label{28.02.1}
\Delta_h^{3} P(y) = 0, \ \ y, h  \in \overline{Y^{(2)}}.
\end{equation}
Hence, $P(y)$ is a continuous polynomial on the group $\overline{Y^{(2)}}$. In particular $P(y)$ is a continuous polynomial on the group $\overline{K^{(2)}}$.
Since   $D$ is a   locally compact  totally disconnected Abelian group,  the group $K$  consists of compact elements (\!\!\cite[(24.17)]{Hewitt-Ross}). Hence, the group  $\overline{K^{(2)}}$ also consists of compact elements. By Lemma \ref{lem4}, we have $P(y)=P(0)=0$ for all  $y\in  \overline{K^{(2)}}$. Since $\psi_j(y)\le 0$ for all $y\in Y$,
it follows from  (\ref{16.09.15.1})   that
\begin{equation}\label{11.02.1}
\psi_1((I+\tilde\alpha)y)=0, \ \ y\in \overline{K^{(2)}}.
\end{equation}

Let us verify that
\begin{equation}\label{08.02.1}
\overline{(I+\tilde\alpha)(K^{(2)})}=\overline{K^{(2)}}
\end{equation}
holds.
Consider the factor-group $X/\mathbb{R}$, and denote by $[x]$ its elements.
Since the subgroup $\mathbb{R}$ is the connected component of zero of the group  $X$, we have
$\alpha(\mathbb{R})=\mathbb{R}$. This implies that $\alpha$
induces a topological automorphism  $\hat\alpha$ on the factor-group $X/\mathbb{R}$ by the formula $\hat\alpha [x]=[\alpha x]$.
Check at the beginning that
\begin{equation}\label{091}
{\rm Ker}(I+\hat\alpha)=\{0\}.
\end{equation}
Let $x_0\in X$ and assume that $[x_0]\in {\rm Ker}(I+\hat\alpha)$. Then
$(I+\hat\alpha)[x_0]=0$,  that is $[(I+\alpha)x_0]=0$.
This implies that  $(I+\alpha)x_0=t_0$
 for some $t_0\in \mathbb{R}$.
It is obvious that
$(I+\alpha)(\mathbb{R})\subset\mathbb{R}$.
Since, in view of
$(\ref{1})$, the restriction of the continuous endomorphism  $I+\alpha$ of the group  $X$ to the subgroup $\mathbb{R}$ is a topological automorphism of the group   $\mathbb{R}$,
we have  $t_0=(I+\alpha) \tilde t$
for some  $\tilde t\in \mathbb{R}$.
Thus,
$(I+\alpha)x_0=(I+\alpha)\tilde t $. In view of $(\ref{1})$, we find from this that $x_0=\tilde t$.
Hence,
  $[ x_0]=0$ and $(\ref{091})$
is proved.

We note that the character group of the factor-group   $X/\mathbb{R}$
in a natural way is topologically isomorphic  with the annihilator $A(Y, \mathbb{R})=K$ (\!\!\cite[(23.25)]{Hewitt-Ross}), and the adjoint of topological automorphism  ${\hat\alpha}$ acts on $K$ as the restriction   of the topological automorphism  $\tilde\alpha$ to
   $K$. For this reason, in view of Lemma \ref{lem4new},
(\ref{091}) implies that the subgroup   $(I+\tilde\alpha)(K)$ is dense in $K$.  Hence,    it follows from this that (\ref{08.02.1}) is fulfilled.

In view of    (\ref{08.02.1}), we find from (\ref{11.02.1}) that $\psi_1(y)=0$
 for all $y\in \overline{K^{(2)}}$. Hence, $\nu_1(y)=1$
 for all $y\in \overline{K^{(2)}}$.

Take an arbitrary element $h_3$ of the group $Y$.
Substitute   $v+h_3$ instead of $v$ in
(\ref{8}). Subtracting
    equation
  (\ref{8}) from the obtaining equation, we get
  \begin{equation}\label{9a}
   \Delta_{h_3}\Delta_{-(I+\tilde\alpha) h_2}\Delta_{-2\tilde\alpha  h_1} Q(v)=0,
\ \ v\in Y.
\end{equation}
Reasoning for the function $Q(y)$ as well as for the function $P(y)$,  we obtain from (\ref{9a}) that the function $Q(y)$ satisfies the equation
\begin{equation}\label{28.02.2}
\Delta_h^{3} Q(y) = 0, \ \ y, h  \in \overline{K^{(2)}}.
\end{equation}
It means that $Q(y)$ is a continuous polynomial on the subgroup  $\overline{K^{(2)}}$, and hence by Lemma \ref{lem4},    $Q(y)=Q(0)=0$ for all $y\in  \overline{K^{(2)}}$.
Then, it follows from  (\ref{16.09.15.1}) that $\psi_2((I+\tilde\alpha)y)=0$
for all $y\in \overline{K^{(2)}}$. Taking into account (\ref{08.02.1}), from this we get that    $\psi_2(y)=0$
 for all $y\in \overline{K^{(2)}}$. Hence, $\nu_2(y)=1$
 for all $y\in \overline{K^{(2)}}$.

So, we have proved that  $\nu_j(y)=1$
for all $y\in \overline{K^{(2)}}$. By Lemma \ref{lem6}, it follows from this that    $\sigma(\nu_j)\subset A(X, \overline{K^{(2)}})$, $j=1, 2$.  It is obvious that $A(X, \overline{K^{(2)}})=\mathbb{R}\times G$. Hence, (\ref{18.04.4}) holds.

2. Let us verify that it suffices to prove the theorem in the case when the group $X$ is of the form $X={\mathbb R}\times G$, where  $G$ is a   locally compact    Abelian group such that all its nonzero elements have order   $2$, and    $\hat\mu_j(y)>0$ for all $y\in Y$.
Since $\sigma(\nu_j)\subset \mathbb{R}\times G$ and $\nu_j = \mu_j* \bar \mu_j$, it is easy to see that we can replace the distributions  $\mu_j$ by their shifts  $\mu'_j=\mu_j*E_{-x_j}$, $x_j\in X$,   in such a way that $\sigma(\mu'_j)\subset \mathbb{R}\times G$, $j=1, 2$. It is clear that  $\nu_j = \mu'_j* \bar \mu'_j$. If we prove that $\nu_j\in \Gamma({\mathbb R})*{\rm M}^1(G)$, it follows from Lemma  \ref{lem5},  applied to the group $\mathbb{R}\times G$, that
$\mu'_j=\gamma_j*\omega_j$, where $\gamma_j \in \Gamma({\mathbb R})$,
 $\omega_j\in {\rm M}^1(G)$, and hence,
$\mu_j=\gamma_j*\omega_j*E_{x_j}$,  $x_j\in X$, $j=1, 2$. Thus, the theorem will be proved.

Denote by $\eta_j$ independent random variables with values in the group $X$ and distributions $\nu_j$, $j=1, 2$. Since the characteristic functions   $\hat\nu_j(y)$ satisfy equation $(\ref{2})$, by Lemma \ref{lem1},  the conditional distribution of the linear form $M_2 = \eta_1 + \alpha\eta_2$ given $M_1 = \eta_1 + \eta_2$ is symmetric. Since $\mathbb{R}$ is the connected component of zero of the group   $X$, and $G=\{x\in D:2x=0\}$, we have
$\alpha(\mathbb{R})=\mathbb{R}$,    $\alpha(G)=G$, and so that
$\alpha(\mathbb{R}\times G)=\mathbb{R}\times G$. Obviously,   the topological automorphism   $\alpha_{\mathbb{R}\times G}$ of the group $\mathbb{R}\times G$ satisfies condition (\ref{1}). We can consider  $\eta_j$ as  independent random variables with values in the group  ${\mathbb R}\times G$. Then  the conditional distribution of the linear form   $M_2 = \eta_1 + \alpha_{\mathbb{R}\times G}\eta_2$ given $M_1 = \eta_1 + \eta_2$ is symmetric. It follows from what has been said that in the proof of the theorem one can assume from the  beginning that $X={\mathbb R}\times G$, where  $G$ is a   locally compact    Abelian group such that all its nonzero elements have order   $2$, and    $\hat\mu_j(y)>0$ for all $y\in Y$.

3. So, let $X={\mathbb R}\times G$, and let $\hat\mu_j(y)>0$ for all $y\in Y$.  It is obvious that $Y=\mathbb{R}\times H$, where the group   $H$  is topologically isomorphic to the character group of the group  $G$. Denote by $x=(t, g)$, where $t\in \mathbb R$, $g\in G$, elements of the group $X$ and by $y=(s, h)$, where $s\in \mathbb R$, $h\in H$, elements of the group $Y$. Since
$\alpha(\mathbb{R})=\mathbb{R}$ and $\alpha(G)=G$, this implies that $\alpha$ acts on the elements of the group $X$ as follows $\alpha(t, g)=(at, \alpha_Gg)$,   where $a\ne 0$, $t\in \mathbb R$, $g\in G$. In view of (\ref{1}),  $a\ne -1$. Note that all nonzero elements of the group $H$ are of order $2$. For this reason $h=-h$ for all $h\in H$. Taking this into account, equation (\ref{2})
   for the characteristic functions $\hat\mu_j(s, h)$   takes the form
\begin{equation}\label{14.04.2}
\hat\mu_1(s_1+s_2, h_1+h_2)\hat\mu_2(s_1+a s_2, h_1+ \tilde\alpha_Gh_2)$$$$=
\hat\mu_1(s_1-s_2, h_1+h_2)\hat\mu_2(s_1-a s_2, h_1+\tilde\alpha_Gh_2), \ \  s_j\in \mathbb{R}, \ \ h_j\in H.
\end{equation}
Setting $h_1=h_2=0$ in (\ref{14.04.2}) we get
\begin{equation}\label{09.02.1}
\hat\mu_1(s_1+s_2, 0)\hat\mu_2(s_1+a s_2, 0)=
\hat\mu_1(s_1-s_2, 0)\hat\mu_2(s_1-a s_2, 0), \ \  s_j\in \mathbb{R}.
\end{equation}
Taking into account the fact that   $\hat\mu_j(s, 0)>0$ for all $s\in \mathbb{R}$, and Lemmas  \ref{lem1} and \ref{lem7}, it follows from (\ref{09.02.1}) that $\hat\mu_j(s, 0)=\exp\{-\sigma_js^2\}$, $s\in \mathbb{R},$  where $\sigma_j\ge 0$,   $j=1, 2$.  Note also that  $\sigma_1+a\sigma_2=0$. Therefore, either   $\sigma_1=\sigma_2=0$ or $\sigma_1>0$ and $\sigma_2>0$.

If $\sigma_1=\sigma_2=0$, then $\hat\mu_1(s, 0)=\hat\mu_2(s, 0)=1$ for all $s\in \mathbb{R}$. Hence, by Lemma  \ref{lem6}, we have  $\sigma(\mu_j)\subset A(X,  \mathbb{R})$, $j=1, 2$. Since $A(X,  \mathbb{R})=G$,  in this case  the theorem is proved. Therefore, we will assume that
 \begin{equation}\label{14.04.3}
\hat\mu_j(s, 0)=\exp\{-\sigma_js^2\}, \ \ s\in \mathbb{R}, \ \ \sigma_j> 0, \ \ j=1, 2.
\end{equation}

By Lemma \ref{lem2}, it follows from (\ref{14.04.3}) that for each fixed $h\in H$ the functions $\hat\mu_j(s, h)$ can be extended to the complex plane $\mathbb{C}$ as entire functions in   $s$. Moreover, equation (\ref{14.04.2}) is valid for all  $s_j\in \mathbb{C}$, $h_j\in H$. Let us check that the functions  $\hat\mu_j(s, h)$ do not vanish in $\mathbb{C}$. Indeed, suppose that  $\hat\mu_1(s_0, h_0)=0$ for some $s_0\in \mathbb{C}$, $h_0\in H$. Note that by Lemma \ref{lem4new},   ${\rm Ker}(I+\alpha_G)=\{0\}$ implies that the subgroup $(I+\tilde\alpha_G)(H)$ is dense in $H$, that is
\begin{equation}\label{11.02.2}
\overline{(I+\tilde\alpha_G)(H)}=H.
\end{equation}
It follows from  (\ref{11.02.2}) that there is a sequence  $h_2^{(n)}\in H$ such that
 \begin{equation}\label{19.02.1}
(I+\tilde\alpha_G)h_2^{(n)}\rightarrow h_0 \ \ \mbox{as} \ \ n\rightarrow\infty.
\end{equation}
Put  $s_1=s_2=s_0/2$, $h_1=\tilde\alpha_Gh_2^{(n)}$, $h_2=h_2^{(n)}$ in (\ref{14.04.2}).
In view of (\ref{14.04.3}), we obtain
\begin{equation}\label{14.04.4}
\hat\mu_1(s_0, (I+\tilde\alpha_G)h_2^{(n)})\exp\{-\sigma_2((1+a)s_0/2)^2\}$$$$=
\hat\mu_1(0, (I+\tilde\alpha_G)h_2^{(n)}) \exp\{-\sigma_2((1-a)s_0/2)^2\}.
\end{equation}
Taking into account that by Lemma \ref{lem2}, the functions $\hat\mu_j(s, h)$ are continuous on $\mathbb{C}\times H$, it follows from (\ref{19.02.1}) that the left-hand side of (\ref{14.04.4})  tends to zero as $n\rightarrow\infty$, whereas  the right-hand side does not.  From the obtained contradiction it follows that the function $\hat\mu_1(s, h)$ does not vanish   in the complex plane $\mathbb{C}$.   For the function $\hat\mu_2(s, h)$ we argue similarly. Assume that
$\hat\mu_2(s_0, h_0)=0$ for some $s_0\in \mathbb{C}$,   $h_0\in H$, and  (\ref{19.02.1}) is satisfied. Then setting  $s_1=s_0/2$, $s_2=s_0/2a$, $h_1=h_2=h_2^{(n)}$  in equation (\ref{14.04.2}),  we get  that the left-hand side of the resulting equality  tends to zero as $n\rightarrow\infty$, whereas  the right-hand side does not. Therefore, the function $\hat\mu_2(s, h)$ also does not vanish in the complex plane $\mathbb{C}$.

It follows from (\ref{14.04.3}) and Lemma \ref{lem2} that for each fixed  $h\in H$,   $\hat\mu_j(s, h)$ are entire functions in    $s$  of the order at most   2. Since the functions $\hat\mu_j(s, h)$ do not vanish, by the Hadamard theorem on the representation of an entire function of finite order, we get the representation
\begin{equation}\label{14.04.5}
\hat\mu_j(s, h)=\exp\{A_{hj}s^2+B_{hj}s+C_{hj}\}, \ \ s\in \mathbb{R}, \ \ h\in H, \ \ j=1, 2,
\end{equation}
where the coefficients $A_{hj}$, $B_{hj}$ and $C_{hj}$ continuously depend on $h$. Since  $\hat\mu_j(s, h)>0$, we conclude that $A_{hj}$, $B_{hj}$, $C_{hj}$ are real. Put
$$
\varphi_j(s, h)=  A_{hj}s^2+B_{hj}s+C_{hj}, \ \ s\in \mathbb{R}, \ \ j=1, 2,
$$
and substitute (\ref{14.04.5}) into (\ref{14.04.2}). It follows from the obtained equation that the functions $\varphi_j(s, h)$ satisfy the equation
$$
\varphi_1(s_1+s_2, h_1+h_2)+\varphi_2(s_1+a s_2, h_1+ \tilde\alpha_Gh_2)$$$$=
\varphi_1(s_1-s_2, h_1+h_2)+\varphi_2(s_1-a s_2, h_1+ \tilde\alpha_Gh_2), \ \  s_j\in \mathbb{R}, \ \ h_j\in H.
$$
Hence,
\begin{equation}\label{14.04.7}
A_{h_1+h_2, 1}(s_1+s_2)^2+B_{h_1+h_2, 1}(s_1+s_2)+C_{h_1+h_2, 1}$$$$+A_{h_1+ \tilde\alpha_Gh_2,2}(s_1+as_2)^2
+B_{h_1+ \tilde\alpha_Gh_2,2}(s_1+as_2)+C_{h_1+ \tilde\alpha_Gh_2,2} $$$$=A_{h_1+h_2, 1}(s_1-s_2)^2+B_{h_1+h_2, 1}(s_1-s_2)+C_{h_1+h_2, 1}$$$$+A_{h_1+ \tilde\alpha_Gh_2,2}(s_1-as_2)^2
+B_{h_1+ \tilde\alpha_Gh_2,2}(s_1-as_2)+C_{h_1+ \tilde\alpha_Gh_2,2}, \ \ s\in \mathbb{R}.
\end{equation}
Equating  the coefficients  of $s_1s_2$ in both  sides  of (\ref{14.04.7}), we get
\begin{equation}\label{14.04.8}
A_{h_1+h_2, 1}+aA_{h_1+ \tilde\alpha_Gh_2,2}=0.
\end{equation}
Note that (\ref{14.04.3}) and (\ref{14.04.5}) imply the equalities
\begin{equation}\label{20.02.1}
A_{01}=-\sigma_1, \ \ A_{02}=-\sigma_2.
\end{equation}
Setting first $h_1=h_2$, then $h_1=\tilde\alpha_Gh_2$ in (\ref{14.04.8})  and taking into account  (\ref{11.02.2}), we conclude from (\ref{20.02.1}) that
\begin{equation}\label{15.04.1}
A_{h1}= -\sigma_1, \ \ A_{h2}= -\sigma_2, \ \ h\in H.
\end{equation}
Equating  the coefficients  of $s_2$ in both sides of (\ref{14.04.7}), we get
\begin{equation}\label{14.04.9}
B_{h_1+h_2, 1}+aB_{h_1+ \tilde\alpha_Gh_2,2}=0.
\end{equation}
It follows from  (\ref{14.04.3}) and (\ref{14.04.5}) that
\begin{equation}\label{20.02.2}
B_{01}=B_{02}=0.
\end{equation}
Setting first $h_1=h_2$, then $h_1=\tilde\alpha_Gh_2$ in (\ref{14.04.9})  and taking into account  (\ref{11.02.2}), we conclude from (\ref{20.02.2}) that
\begin{equation}\label{15.04.2}
B_{h1}=B_{h2}=0, \ \ h\in H.
\end{equation}
Observe that (\ref{14.04.5}) implies that  $\hat\mu_j(0, h)=\exp\{C_{hj}\}$ for all $h\in H$,  $j=1, 2$.
Denote by $\omega_j$ a distribution on the group   $G$ with the characteristic function
\begin{equation}\label{15.04.3}
\hat\omega_j(h)=\exp\{C_{hj}\}, \ \ h\in H, \ \ j=1, 2.
\end{equation}
In view of (\ref{15.04.1}), (\ref{15.04.2}) and (\ref{15.04.3}),   (\ref{14.04.5})  implies the representation
\begin{equation}\label{14.04.11}
\hat\mu_j(s, h)=\exp\{-\sigma_j s^2\}\hat\omega_j(h), \ \ s\in \mathbb{R}, \ \ h\in H, \ \ j=1, 2.
\end{equation}
Denote by $\gamma_j$ a Gaussian distribution on the group   $\mathbb{R}$ with  the characteristic function
$$
\hat\gamma_j(s)=\exp\{-\sigma_js^2\}, \ \ s\in \mathbb{R}, \ \ j=1, 2.
$$
It follows from (\ref{14.04.11}) that $\mu_j=\gamma_j*\omega_j$,  $j=1, 2$.
 \hfill$\Box$
\begin{remark}\label{r1}   Theorem \ref{th1} cannot be strengthened by narrowing the class of distributions which is characterized by the symmetry of the conditional distribution of one linear form given another. Indeed, let   $X={\mathbb R}\times D$, where  $D$ is a locally compact  totally disconnected Abelian group. Set $G=\{x\in D:2x=0\}$.
Let  $\alpha$ be a topological automorphism of the group    $X$ satisfying condition  $(\ref{1})$. Suppose that the restriction of  $\alpha$ to ${\mathbb R}$ is of the form $\alpha (t, 0)=(at, 0)$, where $a\ne 0$, $a\ne -1$. Let $\gamma_j$ be   Gaussian distributions on   ${\mathbb R}$ with the characteristic functions
$$
\hat\gamma_j(s)=\exp\{-\sigma_js^2\}, \ \ s\in  {\mathbb R}, \ \ j=1, 2,
$$  and we will assume that  $\sigma_1+a\sigma_2=0$. Let $\omega_j$, $j=1, 2$, be arbitrary distributions on    $G$. If we consider  $\gamma_j$ and $\omega_j$ as distributions on  $X$ it is easy to see that the characteristic functions    $\hat\gamma_j(y)$ and $\hat\omega_j(y)$ satisfy equation (\ref{2}). Therefore, the characteristic functions  $\hat\gamma_j(y)\hat\omega_j(y)$  also satisfy equation (\ref{2}). Set $\mu_j=\gamma_j*\omega_j$, $j=1, 2$.  Let $\xi_1$ and $\xi_2$ be independent random variables with   values in the group  $X$ and distributions
$\mu_1$ and $\mu_2$.
Since the characteristic functions    $\hat\mu_j(y)$  satisfy equation (\ref{2}), by Lemma \ref{lem1}, the conditional distribution of the linear form  $L_2 = \xi_1 + \alpha\xi_2$ given $L_1 = \xi_1 + \xi_2$ is symmetric.
\end{remark}
\begin{remark}\label{r3}
In   item 3 of the proof of Theorem \ref{th1} we proved the following statement.

{\it Let $X$ be a group of the form   $X=\mathbb{R}\times G$, where  $G$ is a   locally compact    Abelian group such that all its nonzero elements have order   $2$. Let  $\alpha=(a, \alpha_G)$ be a topological automorphism of the group    $X$ satisfying condition $(\ref{1})$. Let $\xi_1$ and $\xi_2$ be independent random variables with   values in   $X$ and distributions
$\mu_j$ such that $\hat\mu_j(y)>0$ for all $y\in Y$.
If the conditional distribution of the linear form    $L_2 = \xi_1 + \alpha\xi_2$ given $L_1 = \xi_1 + \xi_2$ is symmetric, then
$\mu_j=\gamma_j*\omega_j$, where $\gamma_j \in \Gamma({\mathbb R})$,
 $\omega_j\in {\rm M}^1(G)$,
  $j=1, 2$.}

It turns out that  this statement, generally speaking, is not true, if $G$ is an arbitrary   locally compact    Abelian group. Indeed, let $G=\mathbb{R}$, i.e. $X=\mathbb{R}^2$, and let $\alpha(t, g)=(-2t, -2g)$, $(t, g)\in X$, be a topological automorphism of the group    $X$. Since $I+\alpha=-I$,  condition $(\ref{1})$ is satisfied. Put $A_1=\left(\begin{matrix}4&2\\ 2&2\end{matrix}\right)$, $A_2=\left(\begin{matrix}2&1\\ 1&1\end{matrix}\right)$. Denote by $\mu_j$ a distribution on the group $X$ with the characteristic function
\begin{equation}\label{21.04.1}
\hat\mu_j(s, h)=\exp\{-\langle A_j(s, h), (s, h)\rangle\}, \ \ (s, h)\in Y, \ \ j=1, 2.
\end{equation}
Let $\xi_1$ and $\xi_2$ be independent random variables with   values in   $X$ and distributions
$\mu_j$. Substituting (\ref{21.04.1}) into (\ref{2}) it is easy to make sure   that the characteristic functions $
\hat\mu_j(s, h)$ satisfy equation (\ref{2}).
By Lemma \ref{lem1},   the conditional distribution of the linear form
  $L_2 = \xi_1 + \alpha\xi_2$ given $L_1 = \xi_1 + \xi_2$  is symmetric. Obviously, the distributions $\mu_j$ are not represented in the form $\mu_j=\gamma_j*\omega_j$, where $\gamma_j \in \Gamma({\mathbb R})$,
 $\omega_j\in {\rm M}^1(G)$,
  $j=1, 2$.

\end{remark}

The following statement easily follows from the proof of   item 1 of Theorem \ref{th1} and can be considered as an analogue of the Heyde theorem for locally compact  totally disconnected Abelian groups.

\begin{proposition}\label{pr1} Let $X$ be a   locally compact  totally disconnected Abelian group. Set $G=\{x\in X:2x=0\}$. Let  $\alpha$ be a topological automorphism of the group    $X$ satisfying condition $(\ref{1})$. Let $\xi_1$ and $\xi_2$ be independent random variables with   values in   $X$ and distributions
$\mu_1$ and $\mu_2$ with nonvanishing characteristic functions.
If the conditional distribution of the linear form    $L_2 = \xi_1 + \alpha\xi_2$ given $L_1 = \xi_1 + \xi_2$ is symmetric, then
$\mu_j=\omega_j*E_{x_j}$, where
 $\omega_j\in {\rm M}^1(G)$,
$x_j\in X$, $j=1, 2$.
\end{proposition}
\noindent{\it Proof}  We keep the notation used in the proof of Theorem \ref{th1}. Arguing in the same way as in   item 1 of the proof of Theorem \ref{th1}, we arrive at equation  (\ref{9})  from which  we get
\begin{equation}\label{06.03.1}
\Delta_h^{3} P(y) = 0, \ \ y, h  \in \overline{Y^{(2)}}.
\end{equation}
Since   $X$ is a   locally compact  totally disconnected Abelian group,  the group $Y$  consists of compact elements, and the group  $\overline{Y^{(2)}}$ also consists of compact elements. By Lemma \ref{lem4}, we receive from (\ref{06.03.1}) that $P(y)=P(0)=0$  for all  $y\in  \overline{Y^{(2)}}$. Taking into account (\ref{16.09.15.1}), this implies that
\begin{equation}\label{06.03.3}
\psi_1((I+\tilde\alpha)y)=0, \ \ y\in \overline{Y^{(2)}}.
\end{equation}
It is obvious that
\begin{equation}\label{06.03.2}
\overline{(I+\tilde\alpha)(Y^{(2)})}=\overline{Y^{(2)}}.
\end{equation}
In view of (\ref{06.03.2}), we get from (\ref{06.03.3}) that $\psi_1(y)=0$, and hence, $\nu_1(y)=1$ for all $y\in \overline{Y^{(2)}}$. By Lemma \ref{lem6},   this implies that
 $\sigma(\nu_1)\subset A(X, \overline{Y^{(2)}})$. Since $A(X, \overline{Y^{(2)}})=G$, we have $\sigma(\nu_1)\subset G$. For the distribution  $\nu_2$ we reason similarly. The statement of the proposition follows from this. \hfill$\Box$

\section{Main theorem}

Consider a sequence \text{\boldmath $a$}=$(a_0, a_1,a_2,\dots)$, where all $a_j$ are natural and $a_j>1$. Denote by $\Delta_{{\text{\boldmath $a$}}}$   the
group of  ${\text{\boldmath $a$}}$-adic integers. We recall that as a
set $\Delta_{{\text{\boldmath $a$}}}$ coincides
with the Cartesian product $\mathop{\mbox{\rm\bf
P}}\limits_{n=0}^\infty\{0,1,\dots ,a_n-1\}$.
Put $\mathbf{u}=(1, 0,\dots,0,\dots)\in \Delta_{{\text{\boldmath $a$}}}$.
Denote by
 $B$ the subgroup of the group
$\mathbb{R}\times\Delta_{{\text{\boldmath $a$}}}$ of the form
$B=\{(n,n\mathbf{u}): n\in  \mathbb{Z}\}$.
Let
$\Sigma_{{\text{\boldmath $a$}}}$ be the
factor-group
$(\mathbb{R}\times\Delta_{{\text{\boldmath $a$}}})/B$. The group $\Sigma_{{\text{\boldmath $a$}}}$ is called
 the ${\text{\boldmath $a$}}$-{adic solenoid}.   The
group $\Sigma_{{\text{\boldmath $a$}}}$ is  compact, connected and has
dimension 1  (\!\!\cite[(10.12), (10.13),
(24.28)]{Hewitt-Ross}). The character group of the group
$\Sigma_{{\text{\boldmath $a$}}}$ is topologically isomorphic to a discrete
group of the form
$$
 H_{\text{\boldmath $a$}}=
\left\{{m \over a_0a_1 \dots a_n} : \ n = 0, 1,\dots; \ m
\in {\mathbb{Z}} \right\}
$$
(\!\!\cite[(25.3)]{Hewitt-Ross}). To avoid     introducing additional notation we will assume that the character group of the group
$\Sigma_{{\text{\boldmath $a$}}}$ coincides with   $H_{\text{\boldmath $a$}}$. Each topological automorphism
$a$ of the group $\Sigma_\text{\boldmath $a$}$
 is of the following form
$a = f_p f_q^{-1}$ for some mutually prime $p$ and $q$,
where $f_p$ and $f_q$ are also topological automorphisms
 of the group $\Sigma_\text{\boldmath $a$}$. We will identify $a$ with the real number $a={p\over q}$. Since the group $\Sigma_\text{\boldmath $a$}$ is connected, we have  $f_n(\Sigma_\text{\boldmath $a$})=\Sigma_\text{\boldmath $a$}$ for each natural $n$.

It is easy to see that the characteristic function of  a Gaussian distribution $\gamma$ on the group $\Sigma_\text{\boldmath $a$}$ is represented in the form
\begin{equation}\label{e26.2.1}
\hat\gamma(r)=(z, r)\exp\{-\sigma r^2\}, \ \ r\in H_\text{\boldmath $a$},
\end{equation}
where $z\in  \Sigma_\text{\boldmath $a$}$, $\sigma \ge 0$.

Based on Theorem \ref{th1} we will prove the main theorem.
\begin{theorem}\label{th2} Let $X$ be a locally compact    Abelian group with the connected component of zero of dimension $1$.  Set $G=\{x\in X:2x=0\}$. Let  $\alpha$ be a topological automorphism of the group    $X$ satisfying condition $(\ref{1})$. Let $\xi_1$ and $\xi_2$ be independent random variables with   values in   $X$ and distributions
$\mu_1$ and $\mu_2$ with nonvanishing characteristic functions.
If the conditional distribution of the linear form    $L_2 = \xi_1 + \alpha\xi_2$ given $L_1 = \xi_1 + \xi_2$ is symmetric, then
$\mu_j=\gamma_j*\omega_j$, where $\gamma_j\in\Gamma(X)$,
 $\omega_j\in {\rm M}^1(G)$, $j=1, 2$.
\end{theorem}
To prove Theorem \ref{th2} we need the following lemmas.
 \begin{lemma}\label{lem_new1}  {\rm(\!\!\cite{FeTVP1})}    Let  $X$ be a locally compact Abelian group containing no elements of order  $2$, and let  $\alpha$ be a topological automorphism of the group  $X$ satisfying condition $(\ref{1})$. Let $\xi_1$ and $\xi_2$ be
independent random variables with values in   $X$ and distributions
$\mu_1$ and $\mu_2$ with nonvanishing characteristic functions. If the conditional  distribution of the linear form $L_2 = \xi_1 + \alpha\xi_2$ given $L_1 = \xi_1 + \xi_2$ is symmetric, then  $\mu_j\in\Gamma(X)$, $j=1, 2$.
\end{lemma}
\begin{lemma}\label{lem_new2}   Let $X$ be a group of the form
$X = \Sigma_\text{\boldmath $a$}\times G$, where  $G$ is a locally compact Abelian group such that all nonzero elements of $G$ have order
 $2$. Let
$\mu$ be a distribution on $X$ of the form  $\mu =\gamma*\omega$, where $\gamma$ is a symmetric Gaussian distribution on $\Sigma_\text{\boldmath $a$}$ and $\omega \in {\rm M}^1(G)$. Assume that the characteristic function
of the distribution $\omega$ does not vanish. Let $\mu=\mu_1*\mu_2$, where $\mu_j$ are distributions on  $X$. Then   $\mu_j = \gamma_j*\omega_j$, where $\gamma_j \in \Gamma(\Sigma_\text{\boldmath $a$})$, $\omega_j\in {\rm M}^1(G)$, $j=1, 2$.
\end{lemma}
\noindent{\it Proof}  Denote by $x=(z, g)$, where $z\in \Sigma_\text{\boldmath $a$}$, $g\in G$, elements of the group $X$. Then  the group $Y$ can be represented in the form   $Y=H_\text{\boldmath $a$}\times H$, where the group     $H$   is topologically isomorphic to the character group of the group    $G$.
  Denote by $y=(r, h)$, where $r\in H_\text{\boldmath $a$}$, $h\in H$, elements of the group $Y$.
Denote by   $\tau: H_\text{\boldmath $a$}\times H\rightarrow \mathbb{R}\times H$   the embedding $\tau(r, h)=(r, h)$. Set $p=\tilde\tau$. Then $p$ is a continuous homomorphism $p:\mathbb{R}\times G\rightarrow\Sigma_\text{\boldmath $a$}\times G$.   Since $\tau=\tilde p$ and the subgroup $\tau(H_\text{\boldmath $a$}\times H)$ is dense in $\mathbb{R}\times H$, by Lemma \ref{lem4new}, $p$ is a continuous monomorphism. Hence, $p$ induces
a monomorphism from the convolution semigroup  ${\rm M^1}(\mathbb{R}\times G)$ to the convolution semigroup ${\rm M^1}(\Sigma_\text{\boldmath $a$}\times G)$. We will keep the notation $p$ for this monomorphism.
Denote by ${\rm M^1}(p(\mathbb{R}\times G))$ the convolution subsemigroup of ${\rm M^1}(\Sigma_\text{\boldmath $a$}\times G)$ consisting of all distributions of ${\rm M^1}(\Sigma_\text{\boldmath $a$}\times G)$ concentrated on the Borel subgroup $p(\mathbb{R}\times G)$. Then $p$ is
an isomorphism of the convolution semigroups ${\rm M^1}(\mathbb{R}\times G)$ and ${\rm M^1}(p(\mathbb{R}\times G))$.

The characteristic function of the distribution $\mu$ can be written the form
\begin{equation}\label{e25.2.1}
\hat\mu(r, h)=\exp\{-\sigma r^2\}\hat\omega(h), \ \ r\in H_\text{\boldmath $a$}, \ \ h\in H,
\end{equation}
where $\sigma\ge 0$.   Let $M$ be a distribution on the group $\mathbb{R}\times G$ with the characteristic function
\begin{equation}\label{e25.2.2}
\hat M(s, h)=\exp\{-\sigma s^2\}\hat\omega(h), \ \ s\in \mathbb{R}, \ \ h\in H.
\end{equation}
We conclude from (\ref{e25.2.1}) and (\ref{e25.2.2}) that
\begin{equation}\label{1.3.1}
\mu=p(M).
\end{equation}
It follows from this that $\mu\in {\rm M^1}(p(\mathbb{R}\times G))$.
Since $\mu=\mu_1*\mu_2$, we can replace the distributions $\mu_j$ by their shifts $\mu'_j$ in such a way that
\begin{equation}\label{1.3.2}
\mu=\mu'_1*\mu'_2
\end{equation}
and $\mu'_j\in {\rm M^1}(p(\mathbb{R}\times G))$.
We get from (\ref{1.3.1}) and (\ref{1.3.2}) that  $M=p^{-1}(\mu'_1)*p^{-1}(\mu'_2)$. By Lemma  \ref{lem5}, $p^{-1}(\mu'_j)=N_j*\omega_j$, where $N_j \in \Gamma({\mathbb R})$ and $\omega_j\in {\rm M}^1(G)$. This implies that $\mu'_j=p(N_j)*\omega_j$,  $j=1, 2$. To complete the proof, it remains to note that $p(N_j)$ are Gaussian distributions on the group $\Sigma_\text{\boldmath $a$}$. \hfill$\Box$

\bigskip

\noindent{\it Proof of Theorem \ref{th2}}    We divide the  proof of the theorem into several steps.

1. Let us verify that it suffices to prove the theorem in the case when $c_X=\Sigma_\text{\boldmath $a$}$.
By the condition of the theorem the connected component of zero $c_X$ has dimension $1$. It means that we have three possibilities for $c_X$: either $c_X$ is topologically isomorphic to the additive group of real numbers $\mathbb{R}$, or $c_X$ is topologically isomorphic to the circle group $\mathbb{T}$, or  $c_X$ is topologically isomorphic to an \text{\boldmath $a$}-adic solenoid $\Sigma_\text{\boldmath $a$}$. If $c_X$ is topologically isomorphic to $\mathbb{R}$, then $c_X$ is a direct factor of $X$, i.e. the group $X$ is topologically isomorphic to a group of the form  $\mathbb{R}\times D$, where  $D$ is a   locally compact  totally disconnected Abelian group. In this case the statement of the theorem follows from Theorem \ref{th1}. Assume that $c_X$ is topologically isomorphic to $\mathbb{T}$. Note that $\alpha(c_X)=c_X$, and there are only two topological automorphisms   $\alpha$ of the group $\mathbb{T}$, namely  $\alpha=\pm I$. We have in both cases ${\rm Ker}(I+\alpha)\neq\{0\}$. Thus, it is impossible that $c_X$ is topologically isomorphic to   $\mathbb{T}$. It remains to consider the case when $c_X$ is topologically isomorphic to $\Sigma_\text{\boldmath $a$}$.  In order not to complicate the notation, we will assume that $c_X=\Sigma_\text{\boldmath $a$}$.

Note that (\ref{1}) implies that the subgroup $\Sigma_\text{\boldmath $a$}$ contains no element of order 2. Indeed, since $\Sigma_\text{\boldmath $a$}$ is the connected component of zero of the group $X$, we have $\alpha(\Sigma_\text{\boldmath $a$})=\Sigma_\text{\boldmath $a$}$. It is obvious that   $\Sigma_\text{\boldmath $a$}$ can contain no more than one element of order 2. For this reason, if
$z\in \Sigma_\text{\boldmath $a$}$, $z\ne 0$,  and $2z=0$, we have $\alpha z=z$. Hence, $z\in {\rm Ker}(I+\alpha)$ that contradicts (\ref{1}). Thus,
\begin{equation}\label{27.02.1}
\Sigma_\text{\boldmath $a$}\cap G=\{0\}
\end{equation}
is fulfilled.

2.  Put $\nu_j = \mu_j* \bar \mu_j$, and  prove that
\begin{equation}\label{19.04.3}
\sigma(\nu_j)\subset \Sigma_\text{\boldmath $a$}\times G, \ \ j=1, 2.
\end{equation}
Denote the annihilator $A(Y, \Sigma_\text{\boldmath $a$})$ by $K$,  and  first make sure that
\begin{equation}\label{19.04.2}
\sigma(\nu_j)\subset A(X, \overline{K^{(2)}}), \ \ j=1, 2.
\end{equation}
Since $\Sigma_\text{\boldmath $a$}$ is the connected component of zero of the group $X$, the subgroup $K$ consists of compact elements (\!\!\cite[(24.17)]{Hewitt-Ross}). Arguing in the same way as in   item 1 of the proof of Theorem \ref{th1}, we arrive at equation (\ref{8}) from which  equations (\ref{28.02.1}) and (\ref{28.02.2}) follow. We should check that (\ref{08.02.1}) is fulfilled.
Consider the factor-group $X/\Sigma_\text{\boldmath $a$}$ and the induced topological automorphism $\hat\alpha$ on the factor-group $X/\Sigma_\text{\boldmath $a$}$. First we prove that (\ref{091}) holds. Since $\Sigma_\text{\boldmath $a$}$ is a connected group, all endomorphisms $f_n$ are epimorphisms (\!\!\cite[(24.25)]{Hewitt-Ross}). This implies that    the restriction of the continuous endomorphism   $I+\alpha$ of the group $X$  to the subgroup $\Sigma_\text{\boldmath $a$}$ is an epimorphism of the group   $\Sigma_\text{\boldmath $a$}$, and hence, in view of (\ref{1}), this restriction is   a topological automorphism of
$\Sigma_\text{\boldmath $a$}$. Taking this into account, similarly as in the proof of Theorem \ref{th1}, (\ref{091}) follows from (\ref{1}). Using the fact that the character group of the factor-group   $X/\Sigma_\text{\boldmath $a$}$ is topologically isomorphic  with the annihilator $A(Y, \Sigma_\text{\boldmath $a$})=K$ and the adjoint of topological automorphism  ${\hat\alpha}$ acts on $K$ as the restriction   of the topological automorphism  $\tilde\alpha$ to
$K$, we obtain (\ref{08.02.1}) from  (\ref{091}). The final part of the proof of the inclusions (\ref{19.04.2})   is based only on (\ref{08.02.1}) and equations (\ref{28.02.1}) and (\ref{28.02.2}). For this reason we can argue as in the proof of Theorem \ref{th1}.

Let us check now that
\begin{equation}\label{19.04.1}
A(X, \overline{K^{(2)}})=\Sigma_\text{\boldmath $a$}\times G.
\end{equation}
Note that $K^{(2)}=Y^{(2)}\cap K$. Indeed, since $(\Sigma_\text{\boldmath $a$})^{(2)}=\Sigma_\text{\boldmath $a$}$ and $K=A(Y, \Sigma_\text{\boldmath $a$})$, this implies that if $y\in Y$ and $2y\in K$, then $y\in K$. The inclusion $Y^{(2)}\cap K\subset K^{(2)}$ follows from this. The reverse inclusion is obvious. It follows from $K^{(2)}=Y^{(2)}\cap K$ that $A(X, \overline{K^{(2)}})=A(X, \overline{{Y^{(2)}\cap K}})$. Since $\Sigma_\text{\boldmath $a$}$ is a compact subgroup, its annihilator $K$ is an open subgroup of the group $Y$ (\!\!\cite[(23.29)]{Hewitt-Ross}). For this reason $\overline{{Y^{(2)}\cap K}}=\overline{Y^{(2)}}\cap K$. Thus, $A(X, \overline{K^{(2)}})=A(X, \overline{Y^{(2)}}\cap K)$. Observe that the annihilator $A(X, \overline{Y^{(2)}}\cap K)$ coincides with   the closed subgroup generated by $A(X, \overline{Y^{(2)}})$ and $A(X, K)$. We have $A(X, \overline{Y^{(2)}})=G$ and $A(X, K)=\Sigma_\text{\boldmath $a$}$. Hence, $A(X, \overline{Y^{(2)}}\cap K)$ coincides with   the closed subgroup generated by $G$ and $\Sigma_\text{\boldmath $a$}$. Also note that since $G$ is closed and  $\Sigma_\text{\boldmath $a$}$ is compact, the subgroup generated by $G$ and $\Sigma_\text{\boldmath $a$}$ is closed. In view of  (\ref{27.02.1}), this subgroup is $\Sigma_\text{\boldmath $a$}\times G$, and  we conclude that (\ref{19.04.1}) is fulfilled. It follows from (\ref{19.04.2}) and (\ref{19.04.1}) that (\ref{19.04.3}) holds.

3. The next step is to reduce the proof of the theorem to the case when the group $X$ is of the form $X=\Sigma_\text{\boldmath $a$}\times G$, where the ${\text{\boldmath $a$}}$-{adic solenoid}  $\Sigma_{{\text{\boldmath $a$}}}$ contains no elements of order 2,  $G$ is a   locally compact    Abelian group such that all its nonzero elements have order   $2$, and  $\hat\mu_j(y)>0$ for all $y\in Y$. The proof is the same as in    item 2 of the proof of Theorem \ref{th1}. The only difference is we should use (\ref{19.04.3})  instead of (\ref{18.04.4}) and Lemma \ref{lem_new2} instead of Lemma   \ref{lem5}.

4. The final part of the proof is  to derive the statement of the theorem from Theorem \ref{th1}. So, let $X=\Sigma_\text{\boldmath $a$}\times G$,   and  let $\hat\mu_j(y)>0$ for all $y\in Y$.  We have $Y=H_\text{\boldmath $a$}\times H$, where the group   $H$  is topologically isomorphic to the character group of the group  $G$. Denote by $x=(z, g)$, where $z\in \Sigma_\text{\boldmath $a$}$, $g\in G$, elements of the group $X$ and by $y=(r, h)$, where $r\in H_\text{\boldmath $a$}$, $h\in H$, elements of the group $Y$.
 Since the group $\Sigma_\text{\boldmath $a$}$ contains no element of order 2, we have $\alpha(G)=G$. This implies that $\alpha$ acts on the elements of the group $X$ as follows $\alpha(z, g)=(az, \alpha_Gg)$,  $z\in \Sigma_\text{\boldmath $a$}$, $g\in G$, where
 $a$ is a topological automorphism of the group $\Sigma_\text{\boldmath $a$}$. In view of (\ref{1}),  $a\ne -I$. Note that all nonzero elements of the group $H$ are of order $2$. For this reason $h=-h$ for all $h\in H$.
 By Lemma \ref{lem1}, the symmetry of the conditional distribution of the linear form
 $L_2$ given
  $L_1$ implies that the characteristic functions
$\hat\mu_j(y)$ satisfy equation   $(\ref{2})$ which takes the form
\begin{equation}\label{25.02.1}
\hat\mu_1(r_1+r_2, h_1+h_2)\hat\mu_2(r_1+a r_2, h_1+ \tilde\alpha_Gh_2)$$$$=
\hat\mu_1(r_1-r_2, h_1+h_2)\hat\mu_2(r_1-a r_2, h_1+\tilde\alpha_Gh_2), \ \  r_j\in H_\text{\boldmath $a$}, \ \ h_j\in H.
\end{equation}
Setting $h_1=h_2=0$ in (\ref{25.02.1}) we get
\begin{equation}\label{25.02.2}
\hat\mu_1(r_1+r_2, 0)\hat\mu_2(r_1+a r_2, 0)=
\hat\mu_1(r_1-r_2, 0)\hat\mu_2(s_1-a r_2, 0), \ \  r_j\in H_\text{\boldmath $a$}.
\end{equation}
Taking into account the fact that   $\hat\mu_j(r, 0)>0$ for all $r\in H_\text{\boldmath $a$}$, and applying Lemmas  \ref{lem1} and \ref{lem_new1} to the group $\Sigma_\text{\boldmath $a$}$, we find  from (\ref{e26.2.1}) and (\ref{25.02.2}) that
\begin{equation}\label{25.02.6}
\hat\mu_j(r, 0)=\exp\{-\sigma_jr^2\}, \ \ r\in H_\text{\boldmath $a$},
\end{equation}
where $\sigma_j\ge 0$,   $j=1, 2$.  Substituting (\ref{25.02.6}) into (\ref{25.02.2}) we get  $\sigma_1+a\sigma_2=0$. Therefore, either   $\sigma_1=\sigma_2=0$ or $\sigma_1>0$ and $\sigma_2>0$.
If $\sigma_1=\sigma_2=0$, then $\hat\mu_1(r, 0)=\hat\mu_2(r, 0)=1$ for all $r\in H_\text{\boldmath $a$}$. Hence, by Lemma  \ref{lem6}, we have  $\sigma(\mu_j)\subset A(X,  \Sigma_\text{\boldmath $a$})=G$, $j=1, 2$. In this case, the theorem is proved. Thus, we will assume that
 \begin{equation}\label{25.02.3}
\hat\mu_j(r, 0)=\exp\{-\sigma_jr^2\}, \ \ r\in H_\text{\boldmath $a$}, \ \ \sigma_j> 0, \ \ j=1, 2.
\end{equation}

Let $\mu$ be a distribution on a locally compact Abelian group $X$. Then the following inequality
\begin{equation}\label{25.02.4}
|\hat\mu(u)-\hat\mu(v)|^2\le  2(1-{\rm Re}\
\hat\mu(u-v)), \ \ u, v\in Y,
\end{equation}
holds. Embed in a natural way the group $H_\text{\boldmath $a$}\times H$ into the group   $\mathbb{R}\times H$. It follows from  (\ref{25.02.3}) and (\ref{25.02.4}) that the characteristic functions  $\hat\mu_j(r, h)$, $j=1, 2$, are uniformly continuous on the subgroup  $H_\text{\boldmath $a$}\times H$ in the topology induced on  $H_\text{\boldmath $a$}\times H$ by the topology of the group   $\mathbb{R}\times H$. Hence, the characteristic functions $\hat\mu_j(r, h)$ can be extended by continuity from the subgroup   $H_\text{\boldmath $a$}\times H$ to the group $\mathbb{R}\times H$. We retain the notation $\hat\mu_j(s, h)$, where $s\in \mathbb{R}$, $h\in H$, for the continued functions.  It follows from  (\ref{25.02.1}) that the characteristic functions  $\hat\mu_j(s, h)$ satisfy equation (\ref{14.04.2}).   By Theorem \ref{th1}, the representation (\ref{14.04.11}) follows from (\ref{14.04.2}). Hence,
$$
\hat\mu_j(r, h)=\exp\{-\sigma_j r^2\}\hat\omega_j(h), \ \ r\in H_\text{\boldmath $a$}, \ \ h\in H, \ \ j=1, 2.
$$
In view of (\ref{e26.2.1}), the statement of the theorem follows from this.

Reasoning as in Remark \ref{r1}, we make sure that Theorem \ref{th2} cannot be strengthened by narrowing the class of distributions that is characterized by the symmetry of the conditional distribution of one linear form given another. \hfill$\Box$
\begin{remark}\label{r2}  Theorem \ref{th2} is not true if condition (\ref{1}) is not satisfied.   A corresponding example can be constructed  already for the simplest locally compact Abelian group with nontrivial connected component of zero and containing elements of order 2, namely     for $X={\mathbb R}\times \mathbb{Z}(2)$. Denote by $x=(t, g)$, where $t\in \mathbb{R}$, $g\in \mathbb{Z}(2)$,    elements of the group $X$.
The group $Y$ is topologically isomorphic to the group $\mathbb{R}\times \mathbb{Z}(2)$.   Denote by $(s, h)$, where $s\in \mathbb{R}$,
 $h\in \mathbb{Z}(2)$,   elements of the group $Y$. It is obvious that each topological automorphism $\alpha$ of the group   $X$ is of the form
 $\alpha(t, g)=(a t, g)$, where $a\ne 0$.
It is clear that if $a\ne -1$, then ${\rm Ker}(I+\alpha)=\mathbb{Z}(2)$.

The following   statement proved in  \cite[Lemma 4.1]{F_solenoid}.

\noindent{\it   Consider on the group $Y$ the functions
$$
f_j(s, h) = \begin{cases}\exp\{-\sigma_j s^2\}, &s\in \mathbb{R}, \ h=0,\\ \kappa_j\exp\{-\sigma_j' s^2\}, &s\in \mathbb{R}, \   h=1,
\\
\end{cases}
$$
where  $\kappa_j$ are real numbers, $0<\sigma_j'<\sigma_j$ and $0<|\kappa_j|\le\sqrt{\sigma_j'\over \sigma_j}$. Then  $f_j(s, h)$ are characteristic functions of some distributions   $\mu_j$ on the group $X$, $j=1, 2$.}

Let   $\alpha$ be a topological automorphism of the group  $X$   such that $a\ne -1$. Suppose that  numbers  $\sigma_j$ and $\sigma'_j$ are chosen in such a way that
\begin{equation}\label{y6}
\sigma_1+a\sigma_2=0, \ \ \sigma'_1+a\sigma'_2=0.
\end{equation}
Let $\xi_1$ and $\xi_2$ be independent random variables with   values in   the group $X$ and distributions
$\mu_1$ and $\mu_2$. It follows from   (\ref{y6})  that the characteristic functions  $\hat\mu_j(s, h)$ satisfy equation (\ref{2}). Then by Lemma \ref{lem1},
   the conditional distribution of the linear form $L_2 = \xi_1 + \alpha\xi_2$ given $L_1 = \xi_1 + \xi_2$ is symmetric, whereas  $\mu_j\notin \Gamma(\mathbb{R})*{\rm M}^1(\mathbb{Z}(2))$, $j=1, 2$.
\end{remark}

In conclusion, we formulate the following general problem.

\medskip

\noindent{\bf Problem}  {\it Let  $X$ be a   locally compact Abelian group. Set $G=\{x\in X:2x=0\}$. Let  $\alpha$ be a topological automorphism of the group    $X$ satisfying condition $(\ref{1})$. Let $\xi_1$ and $\xi_2$ be independent random variables with   values in   $X$ and distributions
$\mu_1$ and $\mu_2$ with nonvanishing characteristic functions.
Describe  all groups $X$ for which the symmetry of the conditional distribution of the linear form    $L_2 = \xi_1 + \alpha\xi_2$ given $L_1 = \xi_1 + \xi_2$ implies that
$\mu_j=\gamma_j*\omega_j$, where $\gamma_j \in \Gamma(X)$,
 $\omega_j\in {\rm M}^1(G)$,
 $j=1, 2$.}

 \bigskip

 As follows from Proposition  \ref{pr1}, Theorem  \ref{th2} and Lemma \ref{lem_new1} such groups, in particular, are locally compact  totally disconnected Abelian groups,  locally compact Abelian groups with the connected component of zero of dimension $1$, and  locally compact Abelian groups  containing
no elements of order $2$.

\bigskip

\bigskip

\noindent{\bf Conflict of interest} The author states that there is no conflict of interest.

\bigskip

\noindent{\bf Data Availability} Data sharing not applicable to this article as no datasets were generated or analysed during the current study.

\vskip 1 cm

\noindent B. Verkin Institute for Low Temperature Physics and Engineering\\
of the National Academy of Sciences of Ukraine\\
47, Nauky ave, Kharkiv, 61103, Ukraine

\bigskip

\noindent e-mail:    gennadiy\_f@yahoo.co.uk

\end{document}